\documentclass[12pt]{amsart}
\usepackage[all,dvips,arc,curve,color,frame]{xy}
\usepackage{amsmath,amssymb}
\usepackage{url}
\title{Delaunay polytopes derived from the Leech lattice}

\author{Mathieu Dutour Sikiri\'c}
\address{M. Dutour Sikiri\'c, Rudjer Boskovi\'c Institute, Bijenicka 54, 10000 Zagreb, Croatia}
\email{mdsikir@irb.hr}

\author{Konstantin Rybnikov}
\address{K. Rybnikov, Department of Mathematical Sciences, University of Massachusetts at Lowell, Lowell, MA 01854, USA}
\email{Konstantin\_Rybnikov@uml.edu}

\thanks{The authors are thankful for the hospitality of the Hausdorff Research Institute for Mathematics in Bonn, where this research was done. First author has been supported by the Croatian Ministry of Science, Education and Sport under contract 098-0982705-2707.}

\def\QuotS#1#2{\leavevmode\kern-.0em\raise.2ex\hbox{$#1$}\kern-.1em/\kern-.1em\lower.25ex\hbox{$#2$}}

\DeclareMathOperator{\GL}{GL} 
\DeclareMathOperator{\Stab}{Stab}
\DeclareMathOperator{\cov}{cov}
\DeclareMathOperator{\tden}{den}
\DeclareMathOperator{\Iso}{Iso}
\DeclareMathOperator{\conv}{conv}
\DeclareMathOperator{\Min}{Min}
\DeclareMathOperator{\rank}{rank}
\DeclareMathOperator{\perfrank}{perfrank}
\DeclareMathOperator{\Aut}{Aut}
\DeclareMathOperator{\vertt}{vert}

\begin{document}
\newcommand{\KK}{\ensuremath{\mathbb{K}}}
\newcommand{\RR}{\ensuremath{\mathbb{R}}}
\newcommand{\NN}{\ensuremath{\mathbb{N}}}
\newcommand{\QQ}{\ensuremath{\mathbb{Q}}}
\newcommand{\ZZ}{\ensuremath{\mathbb{Z}}}

\newtheorem{theorem}{Theorem}
\newtheorem{proposition}[theorem]{Theorem}
\newtheorem{corollary}[theorem]{Corollary}
\newtheorem{lemma}[theorem]{Lemma}
\newtheorem{problem}[theorem]{Problem}
\newtheorem{conjecture}[theorem]{Conjecture}
\newtheorem{claim}[theorem]{Claim}
\newtheorem{remark}[theorem]{Remark}
\newtheorem{definition}[theorem]{Definition}

\begin{abstract}
Given a lattice $L$ of $\RR^n$, a polytope $D$ is called a
{\em Delaunay polytope} in $L$ if the set of its vertices
is $S\cap L$ where $S$ is a sphere having no lattice points in its interior.
$D$ is called \emph{perfect} if the only ellipsoid in $\RR^n$
that contains $S\cap L$ is exactly $S$.

For a vector $v$ of the Leech lattice $\Lambda_{24}$ we define
$\Lambda_{24}(v)$ to be the lattice of vectors of $\Lambda_{24}$
orthogonal to $v$. We studied Delaunay polytopes of $L=\Lambda_{24}(v)$
for $\|v\|^2 \le 22$.
We found some remarkable examples of Delaunay polytopes in such lattices
and disproved a number of long standing conjectures.
In particular, we discovered:
\begin{enumerate}
\item Perfect Delaunay polytopes of lattice width $4$; previously,
the largest known width was $2$.
\item Perfect Delaunay polytopes in $L$, which can be extended to perfect
Delaunay polytopes in superlattices of $L$ of the same dimension.
\item Polytopes that are perfect Delaunay with respect to two lattices
$L \subset L'$ of the same dimension.
\item Perfect Delaunay polytopes $D$ for $L$ with $|\Aut L|=6|\Aut D|$:
all previously known examples had $|\Aut L|=|\Aut D|$ or $|\Aut L|=2|\Aut D|$.
\item Antisymmetric perfect Delaunay polytopes in $L$, which cannot be
extended to perfect $(n+1)$-dimensional centrally symmetric Delaunay
polytopes.
\item Lattices, which have several orbits of non-isometric perfect Delaunay
polytopes.
\end{enumerate}
Among perfect Delaunay polytopes discovered by us many have large vertex-sets
and sporadic simple groups as their isometry groups.
Finally, we derived an upper bound for the covering
radius of $\Lambda_{24}(v)^{*}$, which generalizes the Smith bound
and we prove that it is met only by $\Lambda_{23}^{*}$, the best known
lattice covering in $\RR^{23}$.
\end{abstract}

\maketitle

\section{Introduction}
 Given an $n$-lattice $L \subset \RR^n$, a
sphere $S=S(c,r)$ of center $c \in \RR^n$ and radius $r$ is called an {\em
empty sphere} for $L$ if there is no $v \in L$ such that $\Vert v-c \Vert
< r$.
A polytope $D=D_L(c)$ in $\RR^n$ (not necessarily of full dimension)
is called a {\em Delaunay polytope} in $L$ if the set of its vertices
$\vertt D$ is $S\cap L$ where $S$ is an {\em empty sphere} for $L$
centered at $c$.
An $n$-dimensional Delaunay polytope $D$ in $L$ is {\em perfect} with
respect to $L$ (or, extreme, as in \cite{DL}) if every linear bijective
transformation $\phi$ of $\RR^n$ that maps $D$ onto a Delaunay polytope
in $\phi(L)$
is a composition of a homothety and an isometry (see \cite{DL,DGL92,integerpaper} for more details on the theory of perfect Delaunay polytopes).
Perfect Delaunay polytopes were first studied by Erdahl in connection with his work on quantum mechanics of many electrons \cite{E75,E92}.

The \emph{norm} of a lattice vector is its squared length.
A vector of $L$ is {\em minimal} if it has the smallest non-zero norm.
Denote by $\Min L$ the set of minimal vectors of $L$.
If $\{v_1, \dots, v_n\}$ is a basis of $L$ and the coordinates of $\Min L$
with respect to $\{v_1, \dots, v_n\}$ determine uniquely $L$ up to
isometries and homotheties, then $L$ is called {\em perfect}.
Perfect Delaunay polytopes are inhomogeneous analogs of perfect lattices.
Perfect lattices have been studied for more than 100 years since their 
introduction (as perfect quadratic forms) by Korkin and Zolotarev \cite{KZ}
in 1873.
The book by Martinet \cite{martinet} contains a wealth of information about
perfect lattices. Not as much is known about perfect Delaunay polytopes.
Up to similarity the unit interval in $\ZZ^1$ and the Gosset
polytope $2_{21}$ in $\mathsf{E}_6$ are the only perfect Delaunay polytopes
in dimension $n\leq 6$ \cite{hyp7}.
We currently know a number of infinite series of perfect Delaunay
polytopes \cite{ErdahlInfinite,Grish,DutourInfinite}, a few sporadic
examples \cite{DGL92} related to highly symmetric lattices in dimensions 12-23,
and a large number of examples in dimensions 7-9 \cite{integerpaper}.
A systematic study of
thousands of $9$-dimensional perfect Delaunay polytopes known  prior to
this paper showed certain uniqueness properties of these polytopes and
their lattices.
All of the known examples satisfied the following conditions:
\begin{enumerate}
\item A polytope can be perfect Delaunay only with respect to one lattice.
\item The vertex set of a perfect Delaunay polytope cannot be a proper subset of the vertex set of another perfect Delaunay polytope of the same dimension.
\item A lattice can have at most one isometry class of perfect Delaunay polytopes.
\item The lattice width of a perfect Delaunay polytope is 2.
\item An antisymmetric perfect Delaunay $n$-polytope uniquely determines a centrally symmetric perfect Delaunay $(n+1)$-polytope (see Section 4)
\item The isometry group of a perfect Delaunay polytope in a lattice $L$ determines the linear isometry group of $L$: the latter is either isomorphic to $\Iso D$ if $D$ is centrally symmetric or is an extension of $\Iso D$ with a group of size $2$ if $D$ is antisymmetric.
\end{enumerate}

We set to find counterexamples for some of these properties and we have found counterexamples for all of them, which often have sporadic simple isometry groups.
We have found them in lattices constructed as sections of the famous Leech lattice $\Lambda_{24}$, which plays a prominent role in geometry, algebra, and number theory (see, almost all chapters of \cite{CS}).
More precisely, given $v\in \Lambda_{24}$ we define 
$\Lambda_{24}(v)$ to be the lattice of vectors of $\Lambda_{24}$ orthogonal to $v$ (see Section 2 for details).
Other well-known lattices that are used in our constructions include the laminated lattices $\Lambda_{n}$ \cite[Chapter 6]{CS} and the lattice $O_{23}$ of O'Connor and Pall \cite{OCP}, also known as the shorter Leech lattice \cite[page 179]{CS}.
As a byproduct of our research, we established the covering radius of $\Lambda^*_{23}$, the lattice holding the covering density record in dimension $23$.

It is instructive to draw a parallel between perfect polytopes and perfect lattices for property (1) and (2) of the above list.
A perfect lattice $L$ may have
the following property: there is a  non-trivial centering $L \subset L'$ such
that $\Min L \subset \Min L'$.
For example, $\mathsf{D}_8\subset \mathsf{E}_8$ with
$\Min \mathsf{D}_8 \varsubsetneq \Min \mathsf{E}_8$
and $\mathsf{A}_9\subset \mathsf{A}_9^2$ with
$\Min \mathsf{A}_9 = \Min \mathsf{A}_9^2$ (see \cite{C,CS}).
In Section \ref{MainDelaunay} we describe a perfect Delaunay
polytope of $\Lambda_{23}=\Lambda_{24}(v_2)$ with $47104$
vertices, which extends to a perfect Delaunay polytope with $94208$ vertices
in the index $2$ superlattice $O_{23}$ of $\Lambda_{23}$.
Similarly, we find a perfect Delaunay polytope with $891$ vertices in $\Lambda_{22}$, which remains a perfect Delaunay in an index $2$ superlattice of $\Lambda_{22}$.
Note that in dimension $n\leq 8$ the set of minimal vectors of a perfect lattice uniquely determines the lattice; furthermore, in dimension $n\leq 7$ such a set cannot be extended to the set of minimal vectors of a denser lattice.

The paper is organized as follows. Sections 1.1 and 1.2 introduce basic notions and terminology. Section 2  introduces the Leech lattice and the lattices $\Lambda_{24}(v)$, computes the covering radius of $\Lambda_{23}^{*}$ and proves an upper bound on the covering radius of $\Lambda_{24}(v)^*$.
Section 3 describes the obtained counterexamples to property (1)-(4), (6).
Section 4 discusses the construction of perfect polytopes by lamination;
finds a counterexample to property (5) and characterizes the cases
for which the construction works.

\subsection{Lattices}

A {\em lattice} $L$ is a subgroup of the vector space $\RR^n$ of the form $\ZZ v_1+\dots+\ZZ v_k$, where  $v_1,\dots, v_k$ are independent vectors. The determinant $\det L$
 is defined as the $k$-dimensional volume of the parallelepiped
 \[\{x_1v_1 +\dots+ x_kv_k\;:\; 0 \le x_i\le 1 \mbox{~for~all~}i\}.\]
For $\KK=\RR$ or $\QQ$ we denote by $L\otimes \KK$ the set $\KK v_1+ \dots + \KK v_k$.
A lattice of rank $k$ is called
a $k$-dimensional lattice or simply a $k$-lattice. The full affine isometry group of a lattice $L$ is denoted by $\Iso L$ and its linear subgroup by $\Aut L$.
If $L$ is a \emph{sublattice} of $L'$, then $L'$ is called a \emph{superlattice} of $L$.
The pair $L \subset L'$, where both lattices are of the same rank, is called a \emph{centering} of $L$.

If $L\subset \RR^n$ is a lattice then the dual lattice $L^*$ is
defined as follows:
\begin{equation*}
L^{*}=\{x\in L\otimes \RR\mbox{~~such~~that~~for~~all~~} y\in L\mbox{~we~have~}\langle
x, y\rangle\in \ZZ\}.
\end{equation*}
A lattice is called \emph{integral} if $\langle x,y\rangle\in \ZZ$ for $x,y\in L$, i.e. if $L\subset L^{*}$.
A lattice is called unimodular if $\det L=1$.
A lattice is \emph{self-dual}, that is $L=L^{*}$, if and only if it is 
integral and unimodular.

If $B^n$ is the $n$-dimensional unit ball centered at the origin and $L$ is an $n$-dimensional
lattice, then the {\em covering radius} $\cov L$ is defined as follows:
\begin{equation*}
\cov L=\min\{\mu \quad : \quad L+\mu B^n \quad \mbox{covers}\quad \RR^n\}.
\end{equation*}
It is easy to see (see, e.g. \cite[Section 2.1.3]{CS}) that $\cov L$ is equal to the maximum circumradius of the Delaunay polytopes of $L$.

A $(n-1)$-dimensional lattice $L'$ of $L$ is called {\em primitive} if $(L'\otimes \RR)\cap L=L'$.
A {\em lamination} of an $n$-dimensional lattice $L$ is a partition of $L$ of
the form $\cup_{k \in \ZZ} (L'+kv)$, where $L'$ is a primitive $(n-1)$-dimensional
sublattice $L'\subset L$ and  $v$ is a  vector in $L \setminus L'$.
Set $\Lambda_0=\ZZ^0$. A laminated $n$-lattice $\Lambda_n$ is defined, up to
similarity, as the densest $n$-lattice that has a lamination
$\cup_{k \in \ZZ} (L'+kv)$, with $L'$ isometric to $\Lambda_{n-1}$
\cite[Chapter 6]{CS}.
In general $\Lambda_n$ is not unique, but it is unique in
dimensions $22$--$24$.
The lattice $\Lambda_{24}$ has $196560$ minimal vectors and is
known as \emph{the Leech lattice}.
The Leech lattice is integral and unimodular and therefore self-dual.
The lattices $\Lambda_{23}$, $\Lambda_{23}^{*}$ have $93150$, respectively $4600$, minimal vectors of norm $4$, respectively $3$.
The $4600$ minimal vectors of $\Lambda_{23}^{*}$ generate an index $2$
sublattice called $O_{23}$.
The lattice $O_{23}$ is known as the \emph{shorter Leech lattice}, as
it was constructed by John Leech from the $23$-dimensional Golay code (after
O'Connor and Pall's work).
The lattice $O_{23}$ is integral and unimodular and therefore self-dual.
 The lattice $\Lambda_{23}$ is an index $2$ sublattice of
$O_{23}$; thus, $\Lambda_{23}$ is an index $4$ sublattice of $\Lambda_{23}^*$.
More precisely, $\Lambda_{23}^* / \Lambda_{23} = \ZZ / 4\ZZ$.

\subsection{Delaunay Polytopes in Lattices}\label{BasicNotions}
Given a lattice $L \subset \RR^n$ a point $x \in \RR^n$
defines a (not necessarily full dimensional)
{\em Delaunay polytope} $D_L(x)$ by
\begin{equation*}
D_L(x) = \conv\left\{v \in L : \|x - v\| = \min_{w\in L}\|x - w\|\right\}.
\end{equation*}
Given a full dimensional Delaunay polytope $D$ denote by $c(D)$ its center
and by $L(D)$ the lattice it affinely generates, i.e. the lattice
generated by the difference between vertices of $D$.
If $D$ is a Delaunay polytope with empty sphere $S(c,r)$ for $L=\ZZ v_1+\dots+\ZZ
v_n$ ($\rank L =n$) then the function
\begin{equation*}
\begin{array}{rcl}
f_D:\ZZ^n         & \rightarrow & \RR\\
x=(x_1,\dots,x_n) & \mapsto     & \Vert \sum_{i=1}^n x_i v_i -c \Vert^2-r^2
\end{array}
\end{equation*}
is a polynomial of degree $2$ such that $f_D(x)\geq 0$ for all $x\in \ZZ^n$ and
$f_D(x)=0$ if and only if $\sum_{i=1}^{n} x_i v_i$ is a vertex of $D$.
The dimension of the cone of quadratic functions
\begin{equation*}
C_D=\left\{\begin{array}{c}
f\;:\; f(x)\geq 0 \mbox{~for all~} x\in \ZZ^n, \mbox{~and~} f(x)=0\\
\mbox{~if and only if~} \sum x_i v_i  \mbox{~is a vertex of~} D
\end{array}\right\}
\end{equation*}
is denoted by $\perfrank D$ and called the \emph{perfection rank} of the Delaunay polytope $D$. $D$ is perfect if and only if its perfection rank is $1$. 
The isometry group of a Delaunay polytope $D=D_L(c)$ is denoted by $\Iso_L D$.
The subgroup $\Aut_{L} D$ is the group of isometries of $L$ preserving $D$.
It can happen that $\Aut_{L} D\not= \Iso_{L} D$ but if $L(D)=L$ then
we have equality.

For a Delaunay polytope $D$ and a $(n-1)$-sublattice $L'$ of $L$,
the non-empty sections of $\vertt D$ by hyperplanes $L'+kv$ are called
{\em laminae}.
The lamination number
(equal to the lattice width number plus $1$, see \cite{JMK})
of $D$ is the minimum
over all primitive $(n-1)$-sublattices $L'$ of the number of laminae.
We found perfect Delaunay polytopes with lamination number $5$, while all the
previously known ones had lamination number $3$.
\cite{JMK} inquired about the possible width of Delaunay polytopes
and conjectured that they cannot have large width.
We expect that there exist Delaunay polytopes with arbitrarily high width.

Given a vector $v\in L\otimes \QQ$, denote by $\tden(v)$ the least common
denominator of its coordinate, i.e. the smallest integer $d>0$ such
that $d v\in L$.
A Delaunay polytope $D$ is either centrally-symmetric with respect to its
circumcenter $c$, or \emph{antisymmetric}, in which case
for any $v \in \vertt D$ we have $2c-v \notin \vertt D$.
Note that $\tden(c(D))=2$ if and only if $D$ is centrally symmetric.
From an antisymmetric perfect Delaunay polytope $D$, it is always
possible to get a centrally symmetric Delaunay polytope by stacking $D$ and
$v-D$ for a suitably chosen vector $v \notin L \otimes \RR$. Often, by varying
the vector $v$ one can ensure that the sphere $S$ around $\vertt D$ and
$v-\vertt D$ contains other points of $L + \ZZ v$, in which case $S \cap (L + \ZZ
v)$ is the vertex set of a centrally symmetric perfect Delaunay polytope.
Many centrally symmetric perfect Delaunay polytopes were
constructed in this way and it was open whether this method
always produces a centrally symmetric perfect Delaunay polytope.
We find an antisymmetric perfect Delaunay polytope for which the construction
does not 
produce  a centrally symmetric perfect Delaunay
polytope and we characterize the cases where it does in Section 4.

We say that a finite, nonempty subset $X$ in $\RR^n$ carries a
\emph{spherical $t$-design} if there is a similarity transformation
mapping $X$ to points on the unit sphere
$S^{n-1} = \{x \in \RR^n : \Vert x \Vert = 1\}$
so that for the spherical measure $d\omega$ on $S^{n-1}$
and for all polynomials $f \in \RR[x_1, \ldots, x_n]$ up to degree $t$
we have
\begin{equation*}
\frac{1}{|Y|}\sum_{y \in Y} f(y) = \frac{1}{\omega(S^{n-1})}\int_{S^{n-1}} f(y) d\omega(y).
\end{equation*}
The maximal possible $t$ is called the {\em strength} of the design
(see, for example, \cite{Ven01} for more details).
A Delaunay polytope define a $0$-design on its empty sphere and it is 
a $1$-design if and only if its circumcenter is equal to its barycenter.
In Section \ref{MainDelaunay} we find many perfect $t$-design associated
to perfect Delaunay polytopes.

The algorithms for Delaunay polytopes used in this study are described
in \cite{DSV08} and reference therein and an implementation is 
available from \cite{Dut09}.
The algorithm for enumerating index $2$ sublattices and
superlattices is exposed in \cite{Ven02}.

\section{The lattices $\Lambda_{24}(v)$}\label{LambdaVlattices}

The \emph{Leech lattice} $\Lambda_{24}$ is a remarkable $24$-lattice,
which can be characterized as the unique $24$-dimensional self-dual
lattice whose non-zero vectors have norm at least $4$. In addition
$\Lambda_{24}$ is \emph{even}, i.e. every $v\in \Lambda_{24}$ has even norm.
The symmetry group $\Aut \Lambda_{24}$ is the Conway group $Co_0$
of order $8315553613086720000$.
Let us say that a vector $v\in \Lambda_{24}$ has type $n$ if it is of norm
$2n$ and has type $n_{a,b}$ if it is also the sum of two vectors of types
$a$ and $b$. $\Aut \Lambda_{24} $ is transitive on the vectors of
following types:
\begin{equation*}
2, 3, 4, 5, 6_{2,2}, 6_{3,2}, 7, 8_{2,2}, 8_{3,2}, 8_{4,2}, 9_{3,3}, 9_{4,2}, 10_{3,3}, 10_{4,2}, 10_{5,2}, 11_{4,3}, 11_{5,2}
\end{equation*}
and this exhausts the list of vectors of norm at most $22$ \cite[Section 10.3.3]{CS}.
By $v_{n}$, respectively $v_{n; a, b}$, we denote a vector of type $n$,
respectively $n_{a,b}$. Note that $2v_2$ is of type $8_{2,2}$.

Let $v$ be a primitive vector of the Leech lattice $\Lambda_{24}$,
i.e. one such that $\frac{v}{d}\notin \Lambda_{24}$ for all natural $d>1$.
Then, we define
\begin{equation*}
\Lambda_{24}(v)=\{x\in \Lambda_{24}\mbox{~such~that~} \langle x, v\rangle=0\}
\end{equation*}
and we have
\begin{equation*}
\QuotS{\Lambda_{24}(v)^*}{\Lambda_{24}(v)}\simeq \QuotS{\ZZ}{\Vert v\Vert^2 \ZZ}.
\end{equation*}
Denote by $\Stab(v)$ the stabilizer of $v\in \Lambda_{24}$
under $\Aut \Lambda_{24}$; the stabilizers of the $17$ vectors $v_n$,
$v_{n;a,b}$ are given in \cite[Table 10.4]{CS} and many of them involve
sporadic simple groups.
Any $d\in \NN$ dividing $\Vert v\Vert^2$ determines uniquely
a lattice $\Lambda_{24}(v, d)$ with
\begin{equation*}
\Lambda_{24}(v)\subset \Lambda_{24}(v,d)\subset \Lambda_{24}(v)^*
\quad {\rm and}\quad 
[\Lambda_{24}(v, d) : \Lambda_{24}(v)] = d.
\end{equation*}
For any such $d$ there is a corresponding integral representation of $\Stab(v)$.
We have $\Lambda_{24}(v,1)=\Lambda_{24}(v)$ and $\Lambda_{24}(v,\Vert v\Vert^2)=\Lambda_{24}(v)^{*}$.

The lattice $\Lambda_{24}(v_2)$
is known as the laminated lattice of dimension $23$, $\Lambda_{23}$.
The lattice $\Lambda_{24}(v_{2}, 2)$, 
denoted by $O_{23}$, is self-dual and is known
as the \emph{shorter Leech lattice}, as
it was constructed by John Leech from the $23$-dimensional Golay code (after
O'Connor and Pall's work).
The stabilizer $\Stab(v_2)$ is the sporadic group $Co_2$.
The group $\ZZ_2\times Co_2$ has three (equivalent over $\QQ$) integral
representations in $\GL_{23}(\ZZ)$ as $\Aut \Lambda_{23}$,
$\Aut \Lambda_{23}^*$ and $\Aut O_{23}$.
Those three representations correspond to the integral
$23$-dimensional representations of $Co_2$
enumerated in \cite{plesken85}.

\begin{theorem}\label{CoveringRadTheo}
Let $v$ be a primitive vector of the Leech lattice and define
\begin{equation*}
r(v)=\sqrt{2-\frac{1}{4\Vert v\Vert^2}}.
\end{equation*}

(i) The covering radius of the lattice $\Lambda_{24}(v)^*$ is at most $r(v)$;

(ii) the only lattice $\Lambda_{24}(v)^*$ with covering radius $r(v)$ is $\Lambda_{23}^{*}=\Lambda_{24}(v_2)^*$.\\
The Delaunay polytopes of maximal circumradius of $\Lambda_{23}^*$ belong to a single orbit of Delaunay polytopes, whose representatives have $64$ vertices and $2688$ symmetries.
\end{theorem}
\proof Denote by $p_v$ the orthogonal projection operator
of $\Lambda_{24}$ onto $\Lambda_{24}(v)\otimes \RR$.
It is proved in the appendix of \cite{martinet} that the dual
lattice $\Lambda_{24}(v)^*$ is equal to the projection $p_v(\Lambda_{24})$
of $\Lambda_{24}$.
We suppose that the covering radius of $\Lambda_{24}(v)^*$ is strictly
greater than $r(v)$, that is that there exists a vector
$w\in \Lambda_{24}(v)\otimes \RR$ such that for every
$x\in \Lambda_{24}(v)^*$ we have
$\Vert x- w\Vert > r(v)$.
Define $v'=\frac{1}{\Vert v\Vert^2} v$ and $w'=w+\frac{1}{2}v'$.
For every $y\in \Lambda_{24}$ set
\begin{equation*}
h_{y}=y-w'=p_v\left(y-w'\right)+\alpha v'=(p_v(y) - w) + \alpha v',
\end{equation*}
where
\begin{equation*}
\alpha=\langle v, h_{y}\rangle=-\frac{1}{2}+\langle v, y\rangle\in \frac{1}{2}+\ZZ.
\end{equation*}
Thus we get
\begin{equation}
\begin{array}{rcl}
\Vert h_{y} \Vert^2
&=& \Vert p_v(y)- w\Vert^2 + \alpha^2 \Vert v'\Vert^2\\
&  >  & r(v)^2 + \alpha^2 \Vert v'\Vert^2\geq 2.
\end{array}
\end{equation}
The inequality $\Vert h_{y}\Vert^2 > 2$ contradicts the fact that the covering radius of $\Lambda_{24}$ is $\sqrt{2}$ (see \cite{parker}, \cite[Chapter 23]{CS}).

Suppose now that $D$ is a Delaunay polytope of $\Lambda_{24}(v)^*$
of center $c$ and circumradius $r(v)$.
Define $c'=c+\frac{1}{2}v'$. By an argument similar to (i) we get that $\Vert y - c'\Vert^2\geq 2$ for every $y\in \Lambda_{24}$.
Thus $c'$ is the center of a Delaunay polytope $D'$ of $\Lambda_{24}$ of circumradius $\sqrt{2}$.
Define $f(x)=\langle x, v\rangle$.
We have $f(c')=\frac{1}{2}$ and $f(v)\in \ZZ$ for $v\in \vertt D'\subset \Lambda_{24}$.

If $v \in \vertt D'$ then we have
\begin{equation*}
\Vert p_v(v)-c\Vert^2=2-\frac{1}{\Vert v\Vert^2}\left(\frac{1}{2}-f(v)\right)^2.
\end{equation*}
Thus in order for $D$ to be a Delaunay polytope of circumradius $r(v)$,
it is necessary that $f(v)=0$ or $1$ for $v\in \vertt D'$.
So $D'$ has lamination number $2$ and the vector $v$ is defined up to
a scalar multiple by the corresponding $2$-lamination.
For a $n$-dimensional polytope $P$ a $2$-lamination in two layers
$L_0, L_1$ corresponds to a partition of $\vertt P$
in two subsets $P_0$ and $P_1$.
If ${\mathcal S}=\{v_1,\dots,v_{n+1}\}$ is a set of $n+1$ independent
vertices of $P$ then the possible partitions $\{P_0,P_1\}$ are
determined by the intersections ${\mathcal S}\cap P_0$.
Thus there are at most $2^{n+1}-2$ $2$-laminations on $P$ and they
can be enumerated by considering all subsets of an independent
set ${\mathcal S}$ of $\vertt P$ and checking if they correspond to
a partition $\{P_0,P_1\}$.
The list of $23$ types $D_{\Lambda_{24}}(c_1)$, \dots, $D_{\Lambda_{24}}(c_{23})$ of Delaunay polytopes of $\Lambda_{24}$ of circumradius $\sqrt{2}$ is
known (see \cite{parker}, \cite[Chapter 23]{CS}).
Given a polytope $D_{\Lambda_{24}}(c_i)$ we enumerate its $2$-laminations;
determine the possible vectors $v$; 
keep the ones such that the projection $p_v(c_i)$ determines
a Delaunay polytope of $\Lambda_{24}(v)^*$ of circumradius $r(v)$.
It turns out that, up to equivalence, only one such vector $v$ satisfies
the required conditions. This vector is $v_2$ and so $\Lambda_{23}^{*}$
is the only lattice meeting the bound. \qed

The lower bound of the above theorem was proved in \cite{Smi88}
for the lattice $\Lambda_{23}^{*}$. The Delaunay polytope of $\Lambda_{24}$
that determines the Delaunay polytope of $\Lambda_{23}^{*}$ of maximal circumradius is named $\mathsf{A}_3^8$ \cite[Chapter 23]{CS}.

The only general method for computing the covering radius of a lattice
is to compute the full Delaunay tessellation.
For $\Lambda_{23}$, respectively $O_{23}$, there are $709$,
respectively $5$, orbits of Delaunay polytopes \cite{DSV08}.
For $\Lambda_{23}^*$ the same program yields several hundred
thousands of orbits before the computation could terminate.

\section{Main Delaunay polytopes of $\Lambda_{24}(v)$}\label{MainDelaunay}

The set $\Min \Lambda_{24}$ consists of $196560$ minimal vectors.
For a given vector $v\in \Lambda_{24}$
and $\alpha\in \ZZ$, define $\Min_{\alpha, v} \Lambda_{24}=\{x\in \Min \Lambda_{24}~:~ \langle x, v\rangle=\alpha\}$.
The following facts are easy to check:
\begin{itemize}
\item[(i)] If $\Min_{0, v} \Lambda_{24}\not= \emptyset$ then $\Min \Lambda_{24}(v)=\Min_{0,v} \Lambda_{24}$.
\item[(ii)] If $\alpha\not= 0$ and $\Min_{\alpha, v} \Lambda_{24}$ is of rank $n$, then it defines a Delaunay polytope of $\Lambda_{24}(v)$ of dimension $n-1$ (see \cite[Lemma 13.2.11]{DL}).
\end{itemize}
We call Delaunay polytopes obtained by this method {\em main}.
If a Delaunay polytope is main then its stabilizer in $\Lambda_{24}(v)$
contains $\Stab(v)$.
Obviously, there is a finite number of main perfect Delaunay polytopes
but we are not able at this point to determine the complete list.
Therefore we limit ourselves to $v$ from the first $17$ types.
We denote by $D(v, N)$ the main Delaunay
polytopes of $\Lambda_{24}(v)$ with $N$ vertices since this notation does
not have ambiguity in the cases considered here.
In Table \ref{tab:PerfectDelaunay} we list the informations about the
full dimensional main perfect Delaunay polytopes in the second column
for the $16$ vector types in first column (The $17$ types, except $8_{2,2}$,
which is covered by type $2$) and their possible extension in
lattices $\Lambda_{24}(v,d)$.
For every main Delaunay $D$ and $d$ for which it admits an extension in
$\Lambda_{24}(v,d)$, we give the number $N$ of vertices,
the denominator $den=\tden(c(D))$ of the circumcenter $c(D)$, the strength $s$
of the spherical $t$-design and the index $ind$ of $L(D)$ in
$\Lambda_{24}(v, d)$ by the symbol ``$d$: $(N, den, s, ind)$''.

The remarkable centrally symmetric perfect Delaunay $D(v_3, 552)$
was first identified in \cite{DGL92}, it defines $276$
equiangular lines \cite{LS},
it is universally optimal \cite{CK} and it gives the facet of maximal
incidence of the contact polytope of $\Lambda_{24}$ \cite{contact}.
It was noted in \cite{DL} that a $22$-dimensional antisymmetric perfect
Delaunay with $275$ vertices is included in $D(v_3, 552)$. This polytope
is $D(v_5,275)$ and the lattice $L_{22}$, which it affinely generates
belongs to the $\QQ$-class of lattices of the irreducible finite subgroup
$(C_2 \times Mc). C_2$ of $\GL_{22}(\ZZ)$ \cite{nebeplesken}.
The polytopes $D(v_5, 275)$ and $D(v_3, 552)$ define spherical $t$-design
of strength $4$, respectively $5$ just like Gosset's $2_{21}$ and $3_{21}$, 
which are perfect Delaunay polytopes in $\mathsf{E}_6$ and $\mathsf{E}_7$.
The set $\Min \Lambda_{23}(v_3)^*$ is equivalent to $D(v_3, 552)$ and
the set $\Min L_{22}^*$ is equivalent to $D(v_5, 275)\cup (2c(D(v_5,275)) - D(v_5, 275))$.
Similarly $3_{21}$ is equivalent to $\Min \mathsf{E}_7^*$
and $\Min \mathsf{E}_6^{*}$ is equivalent to $2_{21}\cup (2c(2_{21})-2_{21})$.

Many lattices $\Lambda_{24}(v)$ have several orbits of perfect Delaunay
polytopes. No such example is known in dimension $n\leq 9$. 
It turns out that for a given vector $v$ of the $17$ cases the strength
of the spherical $t$-design is always the same for all main full
dimensional Delaunay polytopes.
In particular, for vectors of type $2$, $3$, $4$, and $6_{3,2}$
the main Delaunay polytopes $P$ define spherical $t$-designs for
$t=7$, $5$, $3$ and $2$.
Our proof was obtained by direct computation and it would be interesting
to have a less computational proof, for example, using modular forms
in the spirit of the theory of strongly perfect lattices
explained in \cite{Ven01}.

The $22$-di\-men\-sio\-nal Delaunay cell $D(v_{6;2,2},891)$ affinely generates
the $22$-di\-men\-sio\-nal lattice $\Lambda_{22}$.
A remarkable fact is that $|\Aut \Lambda_{22}|=6|\Aut D(v_{6;2,2},891)|$.
For all other known perfect Delaunay polytopes $D$ in a lattice $L$
we have $|\Aut L|=|\Aut D|$ if $D$ is centrally symmetric and 
$|\Aut L|=2|\Aut D|$ if $D$ is antisymmetric.

The perfect main Delaunay polytopes $D(v_3, 552)$, $D(v_3,11178)$ and $D(v_{6;2,2},891)$
admit index $3$, $2$ and $2$ superlattices $\Lambda_{24}(v_3, 3)$,
$\Lambda_{24}(v_3,2)$ and $\Lambda_{24}(v_{6;2,2},2)\cap \Lambda_{22}\otimes \RR$ in which they are still Delaunay polytopes
with the same vertex-set. Note that in \cite{DutRyb} we obtained
Delaunay polytopes with the same property by a different method.
Table \ref{tab:PerfectDelaunay} lists $5$ perfect Delaunay polytopes
$D(v,N)$ that admit superlattices $\Lambda_{24}(v,d)$ in which
these polytopes are proper subsets of perfect Delaunay polytopes.
The first, respectively second, phenomenon is a direct analog of
the relation between the set of minimal vectors of the pair of perfect
lattices $\mathsf{A}_{9}\subset \mathsf{A}^2_9$
and $\mathsf{D}_8 \subset \mathsf{E}_8$.

All the found counterexamples were long suspected to exist.
A technique, which we previously used to no avail
for this task was to get more than $85000$ perfect Delaunay polytopes
in dimension $9$ by using the algorithm of \cite{newalgo}.
It seems that such counterexamples exist only in high dimension.

A Delaunay polytope $D$ in a lattice $L$ is called \emph{basic} if
there are $n+1$ vertices $v_0,v_1,\dots,v_n$ of $D$ such that for every 
vertex $v$ of $D$ there exists $\lambda_i\in \ZZ$ such that
$v=\sum_{i=0}^{n} \lambda_i v_i$ and $1=\sum_{i=0}^{n} \lambda_i$.
A non-basic Delaunay polytope is given in \cite{DG} and a lattice $L$
with $\Min L$ generating $L$ but with no basis of minimal vectors
is given in \cite{CS95}.
We do not know any perfect non-basic Delaunay polytopes; finding
one would be extremely difficult.
The corresponding homogeneous
problem (whether there exists a perfect lattice which is spanned by its
minimal vectors, but where there is no basis of minimal vectors) is also
open.

\begin{table}
\begin{tabular}{|c|c|c|}
\hline
Type of $v$ &\multicolumn{2}{|c|}{Main Delaunay polytopes} \\
\hline
$2$ &\multicolumn{2}{|c|}{$1$: $(47104,4,7,1)$, $2$: $(94208,2,7,1)$}\\
    &\multicolumn{2}{|c|}{$1$: $(4600, 2,7,1)$}\\\hline
$3$ &\multicolumn{2}{|c|}{$1$: $(48600,6,5,1)$}\\
    &\multicolumn{2}{|c|}{$1$: $(11178,3,5,1)$, $2$: $(11178,3,5,2)$}\\
    &\multicolumn{2}{|c|}{$1$: $(552,2,5,1)$, $3$: $(552,2,5,3)$}\\\hline
$4$ &\multicolumn{2}{|c|}{$1$: $(47104,8,3,1)$}\\
    &\multicolumn{2}{|c|}{$1$: $(16192,4,3,1)$, $2$: $(32384,2,3,1)$}\\\hline
$5$ & $1$: $(45100,10,1,1)$\hspace{1.2cm} & $1$: $(19450,5,1,1)$\\\hline
$6_{2,2}$ &\multicolumn{2}{|c|}{$1$: $(22518,6,1,1)$}\\\hline
$6_{3,2}$ &\multicolumn{2}{|c|}{$1$: $(43056,12,2,1)$}\\
         &\multicolumn{2}{|c|}{$1$: $(21528,6,2,1)$}\\
         &\multicolumn{2}{|c|}{$1$: $(6072,4,2,1)$, $2$: $(12144,2,3,1)$}\\\hline
$7$ & $1$: $(41152,14,0,1)$ & $1$: $(22825,7,0,1)$\\
    & $1$: $(7900,14,0,1)$ & \\\hline
$8_{3,2}$  &\multicolumn{2}{|c|}{$1$: $(24576,8,0,1)$, $2$: $(47104,4,7,1)$, $4$: $(94208,2,7,1)$}\\\hline
$8_{4,2}$  &\multicolumn{2}{|c|}{$1$: $(39424,16,1,1)$}\\
          &\multicolumn{2}{|c|}{$1$: $(23608,8,1,1)$}\\
          &\multicolumn{2}{|c|}{$1$: $(9472,16,1,1)$}\\
          &\multicolumn{2}{|c|}{$1$: $(2268,4,1,1)$, $2$: $(4536,2,1,1)$}\\\hline
$9_{3,3}$  & $1$: $(37908,18,1,1)$  & $1$: $(14057,9,1,1)$\\
          & $1$: $(10758,6,1,1)$ & $1$: $(3159,9,1,1)$\\\hline
$9_{4,2}$  &\multicolumn{2}{|c|}{$1$: $(37743,18,0,1)$}\\
          &\multicolumn{2}{|c|}{$1$: $(24035,9,0,1)$}\\
          &\multicolumn{2}{|c|}{$1$: $(10879,6,0,1)$, $3$: $(32384,2,3,1)$}\\\hline
$10_{3,3}$ & $1$: $(25300,10,1,1)$  & $1$: $(4325,5,1,1)$\\\hline
$10_{4,2}$ & $1$: $(25036,10,0,1)$  & $1$: $(3489,5,0,1)$\\\hline
$10_{5,2}$ & $1$: $(36454,20,0,1)$  & $1$: $(24266,10,0,1)$\\
          & $1$: $(11882,20,0,1)$  & $1$: $(3993,5,0,1)$\\\hline
$11_{4,3}$ & $1$: $(35200,22,0,1)$  & $1$: $(24332,11,0,1)$\\
          & $1$: $(12760,22,0,1)$ & $1$: $(4832,11,0,1)$\\\hline
$11_{5,2}$ & $1$: $(34782,22,0,1)$ & $1$: $(24200,11,0,1)$\\
          & $1$: $(13122,22,0,1)$  &\\\hline
\end{tabular}
\label{tab:PerfectDelaunay}
\caption{First column gives the $17$ types of vector of $\Lambda_{24}$
of norm at most $22$ except $8_{2,2}$.
The entries ``$d$: $(N, den, s, ind)$'' in second column correspond
to a main Delaunay polytope $D$ of $\Lambda_{24}(v,d)$ with $N$ vertices,
denominator of circumcenter
$den=\tden(c(D))$, strength $s$ of $t$-design and index $ind$ of $L(D)$
in $\Lambda_{24}(v,d)$}
\end{table}

\subsection{Lamination numbers}
Perfect Delaunay polytopes have lamination number at least $3$
\cite[Theorem 10]{integerpaper} and it is conjectured \cite{Barvinok}
that a $n$-dimensional polytope, whose vertices belong to a lattice $L$
and is free of lattice
point in its interior, has lamination number at most $n+1$.
All known perfect Delaunay polytopes in dimensions $6$, $7$, $8$ have lamination
number $3$.
We will give examples of perfect Delaunay polytopes with lamination number
$5$.
Note that there is no general efficient
method known for determining the lamination number of a given polytope.

\begin{theorem}
The polytopes $D(v_3, 48600)$ and $D(v_3,11178)$ have lamination number $5$.
\end{theorem}
\proof Let us assume that $l(D)\leq 4$ for $D=D(v_3, 48600)$ or $D(v_3,11178)$.
This means that we can find a $22$-dimensional sublattice $L'$ and
four vectors $w_1$, $w_2$, $w_3$, $w_4$ such that
the layers $L_i=w_i+L'$ cover $\vertt D$.

We checked with a computer that $L(D)=\Lambda_{24}(v_3)$,
i.e., the difference vectors of $\vertt D$ generate $\Lambda_{24}(v_3)$.
There exists a linear function $f$ on $\Lambda_{24}(v_3)$ such that
$L'=\ker\,f$ and $f(\Lambda_{24}(v_3))=\ZZ$.
We define an index $2$ sublattice
\begin{equation*}
L'_2=\{w \in \Lambda_{24}(v_3)\mbox{~~such~~that~~}f(w)\in 2\ZZ\}
\end{equation*}
of $\Lambda_{24}(v_3)$ and take $w\in \Lambda_{24}(v_3)$ such that $f(w)=1$.
It is not possible for $L'_2$ or $w+L'_2$ to contain all four layers
$w_i+L'$ since if it were so, then $D$ would not be generating.
So, $L'_2\cap D$ or $(w+L'_2)\cap D$ contains at most $2$ layers.
If one of them contains just one layer, then it is of dimension at most $22$.
By enumerating all index $2$ sublattices of $\Lambda_{24}(v_3)$ we found
that $L'_2\cap D$ and $(w+L'_2)\cap D$ are always $23$-dimensional.
So, $L'_2\cap D$ and $(w+L'_2)\cap D$ are contained in two layers and thus
have lamination number $2$.
We enumerated their $2$-laminations by using the same method as in Theorem
\ref{CoveringRadTheo} and found that each $2$-lamination of $L'_2\cap D$,
respectively $(w+L'_2)\cap D$, induces a lamination
of $(w+L'_2)\cap D$, respectively $L'_2\cap D$
with at least three layers.
So, the lamination number of $D$ is $5$. \qed

\section{Construction of Perfect Delaunay Polytopes by Lamination}\label{CentralConstruction}

In \cite[Lemma 15.3.7]{DL}, \cite{Grish}, \cite{integerpaper}, \cite{EOR}
the following construction of centrally symmetric Delaunay
$(n+1)$-polytopes from Delaunay $n$-polytopes is discussed:
take a Delaunay $n$-polytope $D$ of circumcenter $c$ in
an $n$-di\-men\-sio\-nal lattice $L \subset \RR^{n+1}$ and take its inverted
copy $e_{n+1}-D$ in $L+ e_{n+1}$ for a vector $e_{n+1}\in\RR^{n+1}$.
In order for $\vertt D$ and $e_{n+1} -\vertt D$ to lie on a common sphere
$S$ it is necessary and sufficient that $c - (e_{n+1} -c)$ is orthogonal
to $L\otimes \RR$.
So, up to isometry, $e_{n+1}$ is determined by the
distance between the layer $L$ and the layer $L + e_{n+1}$.
Thus for all $\delta\geq 0$ we define the lattice $L(\delta)=L + \ZZ e_{n+1}$
with $\delta$ the square
Euclidean distance between the layers $L$ and $L+e_{n+1}$ and $2c-e_{n+1}$
orthogonal to $L$.

The theorem below clarifies for which Delaunay polytopes $D$
some other vectors of $L(\delta)$ can lie on $S$ and for which Delaunay
$D$ this cannot happen.

\begin{theorem}\label{ClassifTheo}
Let $D$ be a Delaunay polytope in a $n$-dimensional lattice $L$ of center $c$.
For $i\in \ZZ$, define $D_i=D_L((1-2i) c )$ and denote
by $r_i$ the common distance between $(1-2i)c$ and vertices of $D_i$.
Either:
\begin{itemize}
\item[(i)] For all $i$, $r_i\geq r_0$.
Then $L(0)$ is an index $2$ superlattice of $L$
such that $D'=D_{L(0)}(c)$ is a centrally symmetric Delaunay $n$-polytope
containing $D\cup (2c-D)$ with $\perfrank D'\leq \perfrank D$.
\item[(ii)] Or there exists $i$ such that $r_i < r_0$.
Then there exists $\delta_s>0$ such that $D'=D_{L(\delta_s)} (c')$ with $c'=\frac{1}{2} e_{n+1}$ is a centrally symmetric Delaunay $(n+1)$-polytope containing $D\cup (2c' - D)$ with $\perfrank D'\leq \perfrank D$.
\end{itemize}
\end{theorem}

\proof Define
\begin{equation*}
r_i(\delta)=\sqrt{r_i + \delta \left(i-\frac{1}{2}\right)^2}
\end{equation*}
and $c'=\frac{1}{2}e_{n+1}$.
The sphere circumscribing $D$ and $e_{n+1} - D$
is $S(c', r_0(\delta))$ and we have $r_0(\delta)=r_1(\delta)$.
For $i\in \ZZ$, the set of closest points in layer $L+ i e_{n+1}$ to
$c'$ is 
\begin{equation*}
S_i(c)=i e_{n+1}+\vertt D_L((1-2i) c)
\end{equation*}
and the common distance to $c'$ is $r_i(\delta)$.
If there exists an index $i$ such
that $r_i < r_0$ then there exists $\delta_i>0$ such that
$r_i(\delta_i) = r_0(\delta_i)$ and $S_i(c)$ is outside of 
$S(c', r_{0}(\delta))$ if and only if $\delta \geq \delta_i$.
If one takes $\delta_s=\max_{i\in \ZZ} \delta_i$ then
$S(c', r_{0}(\delta))$ is an empty sphere if and only if $\delta \geq \delta_s$
and $D'=S(c', r_{0}(\delta))\cap L(\delta)$ has more than two layers if and
only if $\delta=\delta_s$.
In that case $L(\delta_s)$ is determined by $L$ and
thus $\perfrank\, D'\leq \perfrank\, D$.

On the other hand, if for all $i$ $r_i\geq r_0$ then 
$\delta_s=0$ and $L(0)$ is actually an $n$-dimensional
superlattice of $L$.
We have $\perfrank\, D'\leq \perfrank\, D$
since $D'$ has more vertices than $D$. \qed

If a Delaunay polytope falls into case (i) then we say that this
Delaunay polytope
is of the {\em first type} and otherwise it is of the {\em second type}.
Given a Delaunay polytope $D$ of the first type we can define for $\delta>0$
a polyhedron 
\begin{equation*}
D_{cyl}=\conv \bigcup_{i\in \ZZ\mbox{~s.t.~} r_i=r_0} (i e_{n+1} + \vertt D_i),
\end{equation*}
whose vertices lie in $L(\delta)$ and belong to a cylinder, empty of lattice
points. It is proved in \cite{E92} that this infinite lattice polyhedron
is in fact 
arithmetically equivalent to the product $\vertt D'\times \ZZ$ where $D'$
is the $n$-dimensional Delaunay polytope of case (i).

\begin{corollary}
Take $D$ a Delaunay polytope of a lattice $L$ of center $c$.

(i) If $\tden(c)=2$ or $4$ then $D$ is of first type.

(ii) If $\tden(c)$ is odd then $D$ is of second type.
\end{corollary}
\proof If $\tden(c)=2$ then $D$ is centrally symmetric and thus of the
first type.
If $\tden(c)=4$, then $D$ is asymmetric and therefore $-D$ is
also a Delaunay polytope. Thus when $\tden(c)=4$ there is a
Delaunay polytope centered at $3c$; $r_i=r_0$ for all $i\in \ZZ$
and by Theorem \ref{ClassifTheo} $D$ is of first type.
If $\tden(c)$ is odd then there exists an index $i$ such
that $(1-2i)c\in L$ and thus $r_i=0$.
So, by Theorem \ref{ClassifTheo} $D$ is of second type. \qed

\smallskip
Let us give two examples, from \cite{DL} of the situation
where $\tden(c(D))$ are odd.
The perfect Delaunay polytopes $2_{21}$ and $D(v_5,275)$ have
$\tden(c(2_{21}))=3$, respectively $\tden(D(v_5, 275))=5$.
Thus they are of second type and the higher dimensional centrally
symmetric Delaunay polytope they define are $3_{21}$ and $D(v_3,552)$.

It was an interesting open question whether there exist
Delaunay polytopes
of first type, which are antisymmetric.
This question was bypassed in \cite[Lemma 15.3.7]{DL} by assuming that 
$D$ is antisymmetric and of second type and it was left open in \cite{Grish}.
The polytopes $D(v_2, 47104)$, $D(v_4, 16192)$ and $D(v_{6;3,2}, 6072)$
have $\tden(c(D))=4$ and thus are of first type and antisymmetric.

Using the method of \cite{newalgo}, we obtained $85000$ perfect Delaunay
polytopes in dimension $9$. 
All the ones of first type were centrally symmetric. All the 
centrally symmetric ones were obtained by the construction of
Theorem \ref{ClassifTheo} but we think that this is not the
case in a large enough dimension.


\begin{thebibliography}{99}




\bibitem{Ven02}
C. Bachoc, B.B. Venkov,
Modular forms, lattices and spherical designs,
in R\'eseaux euclidiens, designs sph\'eriques et formes modulaires,
edited by J. Martinet,
Monographie num\'ero 37 de L'enseignement Math\'ematique, 2001.

\bibitem{BBCGKS}
B. Ballinger, G. Blekherman, H. Cohn, N. Giansiracusa, E. Kelly,
A. Sch\"urmann,
Experimental study of energy-minimizing point configurations on spheres,
to appear in Experiment. Math.

\bibitem{Barvinok}
A. Barvinok, 
A Course in Convexity,
Graduate Studies in Mathematics 54, Amer. Math. Soc. 2002.

\bibitem{CK}
H.~Cohn, A.~Kumar,
Universally optimal distribution of points on spheres,
J. Amer. Math. Soc. 20 (2007) 99--148.

\bibitem{parker}
J.H. Conway, R.A. Parker, N.J.A. Sloane,
The covering radius of the Leech lattice,
Proc. Roy. Soc. London Ser. A 380 (1982) 261--290.

\bibitem{CS}
J.H. Conway, N.J.A. Sloane, Sphere Packings, Lattices and Groups
(third edition), Grundlehren der mathematischen Wissenschaften 290,
Springer--Verlag, 1999.

\bibitem{CS95}
J.H. Conway, N.J.A. Sloane,
A lattice without a basis of minimal vectors,
Mathematika 42 (1995) 175--177

\bibitem{C}
H. S. M. Coxeter,  Extreme forms, Canadian J. Math. 3 (1951) 391--441.

\bibitem{hyp7}
M. Deza, M. Dutour,
The hypermetric cone on seven vertices,
Experiment. Math. 12 (2004) 433--440.

\bibitem{DGL92}
M. Deza, V.P. Grishukhin, M. Laurent,
Extreme hypermetrics and L-polytopes,
in G. Hal\'asz et al. eds, Sets, Graphs and Numbers, Budapest (Hungary), 1991,
60 Colloquia Mathematica Societatis J\'anos Bolyai (1992) 157--209.

\bibitem{DL}
M. Deza, M. Laurent, Geometry of cuts and metrics,
Springer--Verlag, 1997.

\bibitem{DutourInfinite}
M. Dutour, Infinite serie of extreme Delaunay polytopes, 
European J. Combin. 26 (2005) 129--132.


\bibitem{DSV08}
M. Dutour Sikiri\'c, A. Sch\"urmann, F. Vallentin, Complexity
and algorithms for computing Voronoi cells of lattices,
Math. Comp. 78 (2009) 1713--1731.

\bibitem{DG}
M. Dutour Sikiri\'c, V.  Grishukhin,
How to compute the rank of a Delaunay polytope,
European J. Combin. 28 (2007) 762--773.
    

\bibitem{DutRyb}
M. Dutour Sikiri\'c, K. Rybnikov, 
Perfect but not generating Delaunay polytopes,
preprint \url{arXiv:0905.4555}

\bibitem{integerpaper}
M. Dutour Sikiri\'c, R. Erdahl, K. Rybnikov,
Perfect Delaunay polytopes in low dimensions,
Integers 7 (2007) A39.


\bibitem{newalgo}
M. Dutour Sikiri\'c, K. Rybnikov,
A new algorithm in geometry of numbers,
In Proceedings of ISVD-07, the 4-th International Symposium on Voronoi
Diagrams in Science and Engineering, Pontypridd, Wales, July 2007. IEEE
Publishing Services, 2007.

\bibitem{contact}
M. Dutour Sikiri\'c, A. Sch\"urmann, F. Vallentin,
The contact polytope of the Leech lattice,
preprint \url{arXiv:0906.1427}.

\bibitem{E75}
R. Erdahl, A convex set of second-order inhomogeneous polynomials with applications to quantum mechanical many body theory,
Mathematical Preprint \#1975-40, Queen's University, Kingston, Ontario (1975).

\bibitem{E92}
R. Erdahl, A cone of inhomogeneous second-order polynomials,
Discrete Comput. Geom. 8 (1992) 387--416.

\bibitem{EOR}
R.M. Erdahl, A. Ordine, K. Rybnikov,
Perfect Delaunay Polytopes and Perfect Quadratic Functions on Lattices, 
Integer points in polyhedra---geometry, number theory, representation theory, algebra, optimization, statistics,  
Contemporary Mathematics 452, American Mathematical Society, (2008) 93--114.



\bibitem{ErdahlInfinite}
R. Erdahl, K. Rybnikov, An infinite series of perfect quadratic forms and big Delaunay simplices in $\ZZ^n$,
Tr. Mat. Inst. Steklova 239 (2002), Diskret. Geom. i Geom. Chisel, 170--178; translation in Proc. Steklov Inst. Math. 239 (2002) 159--167.

\bibitem{Grish}
V. Grishukhin, Infinite series of extreme Delaunay polytopes,
European J. Combin. 27 (2006) 481--495.




\bibitem{JMK}
J.-M. Kantor, Lattice polytope: some open problems, AMS Snowbird Proceedings.




\bibitem{KZ}
A. Korkine, G. Zolotareff,
Sur les formes quadratiques,
Math. Ann. 6 (1873) 366--389.


\bibitem{LS}
P.W.~Lemmens, J.J.~Seidel, Equiangular lines,
J. Algebra 24 (1973) 494--512.


\bibitem{martinet}
J. Martinet, Perfect lattices in Euclidean spaces, Springer, 2003.

\bibitem{nebeplesken}
G. Nebe, W. Plesken, Finite rational matrix groups,
Mem. Amer. Math. Soc. 1995.


\bibitem{OCP}
R.E. O'Connor, G. Pall,
The construction of integral quadratic forms of determinant $1$,
Duke Math. J. 11 (1944) 319--331.

\bibitem{plesken85}
W. Plesken,
Finite unimodular groups of prime degree and circulants,
J. of Algebra 97 (1985) 286--312.


\bibitem{Smi88}
W. Smith,
PhD thesis: studies in computational geometry motivated by mesh generation,
Department of Applied Mathematics, Princeton University, 1988.

\bibitem{Ven01}
B.B. Venkov, R\'eseaux et designs sph\'eriques,
in R\'eseaux euclidiens, designs sph\'eriques et formes modulaires,
edited by J. Martinet,
Monographie num\'ero 37 de L'enseignement Math\'ematique, 2001.


\smallskip

\centerline{\sc Software}

\bibitem{Dut09}
M. Dutour Sikiri\'c, polyhedral, \url{http://www.liga.ens.fr/~dutour/Polyhedral/}.


\end{thebibliography}
\end{document}